\input amstex
\def\iso{\buildrel\sim\over\rightarrow}
\documentstyle {amsppt}
\document
\nologo
\magnification=1200
\NoBlackBoxes
\pageheight{18cm}


\bigskip
\centerline{\bf COMPOSITION OF POINTS}

\medskip

\centerline{\bf AND MORDELL--WEIL PROBLEM FOR CUBIC SURFACES}

\medskip

\centerline{\bf D. Kanevsky${}^{1,2}$, Yu. Manin${}^2$}

\medskip
\centerline{\it ${}^1$T.J. Watson Research Center, P.O. Box 218, 23-116A,}

\centerline{\it Yorktown Heights, New York 10598, US}
\centerline{\it ${}^2$Max--Planck--Institut f\"ur Mathematik, Bonn, Germany}

\bigskip

{\bf Abstract.} Let $V$ be a plane smooth cubic curve
over a finitely generated field $k.$ The Mordell--Weil theorem
for $V$ states that
there is a finite subset $P\subset V(k)$ such that the whole $V(k)$
can be obtained from $P$ by drawing secants and tangents
through pairs of previously constructed points and consecutively
adding their new intersection points with $V.$ Equivalently,
the group of birational transformations of $V$ generated
by reflections with respect to $k$--points is finitely generated. 
In this paper, elaborating an idea from [M3], we establish a Mordell--Weil type 
finite generation result for some birationally trivial cubic surfaces $W$.
To the contrary, we prove that the birational automorphism
group generated by reflections cannot be finitely generated if $W(k)$
is infinite.

\bigskip

\centerline{\bf \S 1. Introduction}

\medskip

{\bf 1.1. Composition of points.} Let $V$ be a cubic hypersurface without
multiple components over a field $k$ in $\bold{P}^d,\ d\ge 2.$
Three points $x,y,z\in V(k)$ (possibly coinciding)
are called {\it collinear} if either $x+y+z$ is the intersection
cycle of $V$ with a line in $\bold{P}^d$ (with correct multiplicities),
or $x, y, z$ lie on a $k$--line belonging to $V$.
If $x,y,z$ are collinear, 
we write $x=y\circ z.$ Thus $\circ$ is a (partial and multivalued)
composition law on $V(k).$ We will also consider its restriction
on subsets of $V(k),$ e.g. that of smooth points.

\smallskip

If $x\in V(k)$ is smooth, and does not lie on a hyperplane component
of $V$, the birational map $t_x:\ V\to V,\ y\mapsto x\circ y,$
is well defined. It is called reflection
with respect to $x$. Denote by $\roman{Bir}\ V$ the full group
of birational automorphisms of $V.$

\smallskip

The following two results summarize the properties of $\{ t_x\}$
for curves and surfaces respectively. The first one is classical,
and the second is proved in [M1], Chapter V.

\smallskip

\proclaim{\quad 1.2. Theorem} Let $V$ be a smooth cubic curve. Then:

(a) $\roman{Bir}\ V$ is a semidirect product of a finite group
and the subgroup consisting of products of an even
number of reflections $\{ t_x\ |\ x\in V(k)\}.$

\smallskip

(b) We have identically
$$
t_x^2=(t_xt_yt_z)^2=1 \eqno{(1.1)}
$$
for all $x,y,z\in V(k).$

\smallskip

If in addition $k$ is finitely generated over a prime field, then:

\smallskip

(c) $\roman{Bir}\ V$ is finitely generated.
\smallskip

(d) All points of $V(k)$ can be obtained from a finite subset of them
by drawing secants and tangents and adding the intersection
points.

\endproclaim

\smallskip

\proclaim{\quad 1.3. Theorem} Let $V$ be a minimal smooth cubic surface
over a perfect non--closed field $k.$ Then:

(a) $\roman{Bir}\ V$ is a semi--direct product of the group of
projective automorphisms and the subgroup generated by
$$
\{ t_x\ |\ x\in V(k)\}\roman{\ and\ } \{ s_{u,v}\ |\ u,v\in V(K);\
[K:k]=2;\
u,v\ \roman{are\ conjugate\ over\ }k\}
$$
where
$$
s_{u,v}:=t_ut_{u\circ v}t_v,
$$
and $u,v$ do not lie on lines of $V.$

\smallskip

(b) We have identically
$$
t_x^2=(t_xt_{x\circ y}t_y)^2=(s_{u,v})^2=1,\ st_xs^{-1}=t_{s(x)},
\eqno{(1.2)}
$$
for all pairs $u,v$ not lying on lines in $V$, and projective
automorphisms $s.$

\smallskip

(c) The relations (1.2) form a presentation of $\roman{Bir}\ V.$

\endproclaim

\smallskip

We remind that $V$ is called {\it minimal} if one cannot blow down
some lines of $V$ by a birational morphism defined over $k.$ The
opposite class consists of {\it split} surfaces upon which all
lines are $k$--rational.

\medskip

{\bf 1.4. Main results of the paper.} Although Theorems 1.2
and 1.3 look very similar, 
there is an important difference between finiteness
properties in one-- and two--dimensional cases.

\smallskip

Basically, (1.1) means only that $x+y:=e\circ (x\circ y)$ is an Abelian
group law with identity $e$: see [M1],
Theorem I.2.1. The statements c) and d) of
the Theorem 1.2 additionally assert that this group is finitely
generated. Therefore, (1.1) generally is not a complete
system of relations between $\{ t_x\}$.

\smallskip

On the contrary, (1.2) is complete, and in \S 2 we will see that
this prevents $\roman{Bir}\,V$ from being
finitely generated if $V(k)$ is infinite. This answers one of the questions
raised in [M3].

\smallskip

Therefore, any reasonable analog of the Mordell--Weil
problem must address the problem of finite generation
for $(V(k),\circ )$ or of
quotients of $V(k)$ with respect to
various equivalence relations compatible with $\circ$.
This is the subject of \S\S 3--5.

\smallskip

As in [M1], Chapter II, we can start with the universal
equivalence relation $U$. By definition, this
is the finest equivalence relation compatible
with collinearity and such that $\circ$ induces
a well defined operation on $V(k)/U$ also denoted $\circ$. Then one
of the Mordell--Weil type questions asks about
finite generation (= finiteness) of the 
CH--quasigroup $(V(k)/U,\circ )$ (see [M1], Chapter I.).

\smallskip

In \S 3 and \S 4 we give a  description of $U$
refining earlier results of [M1]. Consider
the set of intersections of $V$ with tangent planes
at points of $V(k)$
and add to it all images of these curves with respect to the group
generated by all $t_x, x\in V(k)$. Then one class of $U$ consists
of points that can be pairwise joined by a chain of curves
belonging to this set of curves. This is the content of
Theorem 3.3 below. We then discuss various versions of finite generation
of $(V(k),\circ )$. One essential choice is whether to allow to
apply $\circ$ only to the different previously
constructed points (for minimal surfaces, the result
will then be uniquely defined). Another option
giving more flexibility is to allow expressions
$x\circ x$ and treat them as multivalued, thus adding at one step
all the intersection points of $V$ with a 
tangent plane at $x$. Finally, in \S 4 we
extend the group--theoretic description of $U$
given in [M1], II.13.10.

\smallskip

The results of \S 3 and \S 4 are essentially
algebraic and do not add any new cases of finite
generation of $(V(k),\circ )$ to the short list of locally compact
local fields already treated in [M1]. (In fact,
[M1] proves the finiteness of $V(k)/U$ over such fields
by establishing that $V(k)$ is covered by a finite number of 
sets of the form $(x\circ x) \circ (y\circ y)$).

\smallskip

In \S 5 we study modified composition laws of points 
introduced in [Ma3]. The idea behind this development
is  to reinterpret the classical theorem on the structure
of abstract projective planes as a finiteness
result.

\smallskip

Namely, let $k$ be a finitely generated field.
Start with a finite subset $S\subset \bold{P}^2(k)$
and add to it pairwise intersections of all
lines passing through two points of $S$
thus getting a new finite set $S^{\prime}.$
Apply the same procedure to $S^{\prime}$, and so on.
If $S$ is large enough, in the limit we will
get the whole $\bold{P}^2(k)$. This easily follows
from the fact that if we start
with $S$ consisting of $\ge 4$ points in general position, 
in the limit we will get
an abstract projective plane satisfying
the Desargues axiom and therefore coinciding
with $\bold{P}^2(k^{\prime})$ for $k^{\prime}\subset k$
up to a projective coordinate change.

\smallskip

A trick, first introduced in [M3], allowed us to
translate this remark into a finiteness theorem
for $V(k)$ assuming the existence of a birational
morphism $p:\,V\to \bold{P}^2$ defined over $k$.
However, this required dealing with modified composition
laws: roughly speaking, instead of looking at the collinearity
relation induced by that in $\bold{P}^3$, we now have
to use the collinearity relations determined by
the morphism $p$.

\smallskip

In this paper we make some steps towards eliminating
this complication. Although the final result
falls short of what we would like to prove,
we feel that the connection and analogies with the theory
of abstract projective planes deserve further
study.

\medskip
 
{\bf Acknowledgement.} The first named author would like to thank V. Berkovich and J-L. Colliot-Th\'el\`ene for useful discussions. The work was partially supported by the Humboldt
Foundation during the author's stay at the Max--Planck--Institut f\"ur
Mathematik.

\bigskip

\centerline{\bf \S 2. Cardinality of generators of subgroups in a
reflection group}

\medskip

{\bf 2.1. Notation.} We shall call an {\it abstract cubic } a set $S$ with a ternary relation $L\subset S \times S \times S $, 
satisfying the following axioms:

\smallskip

  (a) $L$ is invariant with respect to permutations of factors.

\smallskip

  (b) If $(x,y,z) , (x,y,z') \in L$ and $x \ne y$, then $z = z^\prime$.

\smallskip

The {\it reflection group} $G_S$ of an abstract cubic $S$ is generated by symbols
$t_x, x \in S$, subject to the following relations:
\smallskip
$t_x^2 = 1$ for all $x \in S$;
\smallskip
$(t_xt_yt_z)^2 = 1$ for all $(x,y,z) \in L$.

\medskip

The following result is proved in  [K1].
\proclaim{\quad 2.2. Theorem} (a). Any element of finite order in $G_S$ is conjugate to either $t_x$ or
to $t_xt_yt_z$ for appropriate $x,y,z \in S$.
\smallskip
Let $S$ be given effectively and $L \subset S \times S \times S $ be
decidable. Then:
\smallskip
(b) The word problem in $G_S$ is decidable.
\smallskip
(c) The conjugacy problem in $G_S$ is decidable.
\endproclaim
\smallskip

The proof is based on a direct description of $G_S$ as a limit of
amalgamated sums.
In [K2] it is shown that $S$ can be sometimes reconstructed from $G_S$. Moreover, under some additional assumptions it is proved
that $\roman{Aut}\,G_S$ is generated by $G_S$ and permutations 
of $S$ preserving $L$.
\smallskip
A different interesting description of $G_S$ and another proof of the Theorem 2.2 is given in [P].

\smallskip

For the purposes of our paper we need the following  description of $G_S$ that
is a special case of the general structure theorem 1.4 in [K1].

\medskip

\proclaim{\quad 2.3. Structural Theorem}
Let $x \in S$ be an arbitrary point and $S' := \{S \setminus {x}\}$. Then
$ G_S$ is canonically isomorphic to $ G_{S'}*_{\Pi}K$ (the free product of
$G_{S'}$ and $K$ with the amalgamated subgroup $\Pi$).
The groups in this product can be described as follows.
\smallskip
(a) $G_{S^{\prime}}$ is the reflection group of the cubic $S^{\prime}$
with the ternary relation induced by $L$ on $S^{\prime}$.
\smallskip
(b) The amalgamated subgroup $\Pi$ is a free group generated by free
generators $a_{u,v} = t_ut_v$ for all distinct pairs $u,v \in S^{\prime}$
such that $(u,v,x) \in L$ and $u < v$ (for some fixed ordering of $S$).
\smallskip
(c) $K \iso Z_2*Z_2*\dots*Z_2$.  Generators of the subgroups $Z_2$ in this free
product are $t_x$ and $t_xa_{u,v}$.
\smallskip
(d) $\Pi$ is of index 2 in  K. The quotient group $K/\Pi$ is generated by the
class $t_x$.
\endproclaim

This structural result leads to the following auxiliary statement, which we
need to prove our main results in this section.

\medskip
\proclaim{\quad 2.4. Definition--Lemma}
(a) In the situation of Theorem 2.3 a family  $W =\langle R_1 , t_x, R_2, t_x,
\dots , t_x , R_n \rangle$, where $R_i \in G_{S^{\prime}}$, is called a reduced
$t_x$--partition of $g = R_1t_xR_2t_x \dots R_n$ if $R_i \notin \Pi $ for $1 <
i < n $.
\smallskip
(b) Let $W$ be a reduced $t_x$--partition of $g \in G_S $. Let us define
$ord_x (g)$ as the number of $t_x$ in $W$. This number depends on $g$ and $x$
and is the same for different reduced $t_x$--partitions of $g$.
\smallskip
(c) Let $g \in G_S$ be such that $ord_x(g)=0$. Then $g \in G_{S^{\prime}}$.
\smallskip
(d) $ord_x(g_1g_2) \equiv (ord_x(g_1) + ord_x(g_2))\,\roman{mod}\,2$.
\smallskip
(e) $ord_x(aga^{-1}) \equiv ord_x(g)\,\roman{mod}\,2$ for any $a,g \in G_S$.
\smallskip
(f) Let $g \in G_S$. We put $\delta(g) := \{x \in S\, |\, ord_x(g) \neq 0 \}$.
The set $\delta(g)$ is finite.\smallskip
(g) $\delta(g_1g_2) \subset \delta(g_1) \cup \delta(g_2)$.
\smallskip
(i) Let $\langle h_1, h_2, \dots\rangle$ be a family generating a subgroup
$H$. Then $\cup_{h \in
H}\delta(h) = \cup_i \delta(h_i)$ .
\endproclaim

\smallskip

Now we can formulate the main theorem of this section.

\medskip

\proclaim{\quad 2.5. Theorem}
Consider a subgroup $H = \langle g_1, g_2, \dots g_n, \dots\rangle\subset G_S$  generated by
an infinite family of elements such that $\delta = \cup_{i}
\delta(g_i)$ is infinite. Then $H$ is not finitely generated.
\endproclaim

\smallskip

{\bf Proof.} Assume that $H$ is finitely
generated by $h_1, \dots ,h_k$. Then
$\delta^{\prime} = \cup_{i=1,..k} \delta(h_i)$ is finite.  Therefore there
exists some $g_r$ and $x \in S$ such  that the following holds: $ord_x(g_r)
\neq 0$ and $x \notin \delta^{\prime}$. Hence $ord_x(h_i)=0$ for all
$i=1,\dots ,k$. By 2.4(i),  $H \subset G_{S^{\prime}}$ if $H$ is generated by
$h_1, \dots ,h_k$. Since $ord_x(g_r) \neq 0$, $g_r \notin H$. This
contradiction proves the theorem.

\smallskip

The following extension of Theorem 2.5 can be applied to 
various subgroups of
$\roman{Bir}\,(V)$.

\medskip
\proclaim{\quad 2.6. Corollary}
Let $G_S$ be the reflection group of an abstract cubic $S$ and let $W$ be
the group of permutations of $S$, preserving its ternary relation $L$.
Let $G \cong W*G_S$ be the semi-direct product of  $W$ and $G_S$, such that
$wt_xw^{-1} = t_{w(x)}$ for any $w \in W$ and $x \in S$.
Let a subgroup $H \subset G$ be generated by a finite subgroup $W^{\prime}
\subset W$ and an infinite family of elements $g_i \in G_S$ such that
$\delta = \cup_{i=1,2,..} \delta(g_i)$ is infinite. Then $H$ is not
finitely generated.
\endproclaim

\smallskip

{\bf Proof.}
Let us assume that $H$ is generated by a finite number of elements $h_1,
\dots ,h_k \in G_S$ and a finite number of elements from $ W^{\prime}$. Let
$\delta = \cup_{i}\delta(g_i)$ and let $\delta^{\prime} = \cup_{w
\in W^{\prime}} w\delta$ be obtained by applications of all $w \in
W^{\prime}$ to $\delta$. Since $\delta$ and $W^{\prime}$ are finite,
$\delta^{\prime}$ is finite. Therefore there exists a generator $g_i \in H$
such that $\delta(g_i) \nsubseteq \delta^{\prime}$. Therefore, $g_i \notin
G_{\delta^{\prime}}$, i.e. it cannot be obtained as a product of elements
from ${h_1,\dots , h_k}$ and $w \in W'$. This contradiction proves the
corollary.
\medskip
{\bf 2.7. Examples.}
In the situation of Theorem 1.3 assume that $S=V(k)$ is infinite.
Then  the following subgroups of $\roman{Bir}\,(V)$ cannot be finitely generated:
\smallskip
(a) $\roman{Bir}\,(V)$, $B(V) := \langle t_x\, |\, x \in V(k)\rangle$ and 
$$
G := \langle t_x , s_{u,v} \,|\, x \in
V(k),\, u,v \in V(K);\, [K:k] = 2; u,v\ \roman{are\ conjugate\ over}\ k\rangle .
$$
(b) The commutant of any of subgroups described in (a).
\smallskip
(c) Let $B_0(V)$ denote the normal subgroup of $B(V)$ generated by elements
of the form $t_xt_yt_zt_{x'}t_{y}t_{z'}$, where $(x,y,z)$ and
$(x^{\prime}, y, z^{\prime})$ run through  triples of
collinear points of $V(k)$. This subgroup was introduced in  [M1] (II.13.9;
beware of a misprint there: the second $y$ carries a superfluous prime).
It is closely related to the universal equivalence on $V(k)$ (see
the section 4 below).
\smallskip
(d)  Let $B_1(V)$ denote the normal subgroup of $B(V)$ generated by
elements of the form $t_xt_yt_z$, where $(x,y,z)$ run through all possible
triples of collinear points of $V(k)$. This subgroup was introduced in 
[M2]  because it is closely related to some admissible equivalence
relations on
$V(k)$
\smallskip
We will now show that Theorem 2.3 implies the statement  for the
case (b). Other cases will be discussed later and stronger statements
will be proved.
\smallskip
{\bf Proof of (b)}.
The commutant of $B(V)$ contains elements
$t_xt_yt_x^{-1}t_y^{-1}=t_xt_yt_xt_y=a_{x,y}^2$. Let us consider an
infinite family of elements $a_{{x_i}{y_i}}^2, x_i, y_i \in S$. The statement
for (b) will follow if we show that $ \delta(a_{x,y}^2)$ contains $x$, since it
will follow that $\cup_i \delta(a_{{x_i}{y_i}})$ is infinite. For this it is enough to show that $t_xt_yt_xt_y$ has the
following reduced $t_x$--partition: $(t_x, t_y, t_x, t_y)$. Indeed, in the
notation of 2.4 one has to check that $t_y \notin \Pi $. But this
fact follows immediately from 2.3 if one notes that $\Pi$ is a free group
(hence it contains no nontrivial elements of finite order) and $t_y^2 = 1$. This
implies that $ord_x (a_{x,y}^2) > 0$, i.e. $x \in  \delta(a_{x,y}^2)$.
Q.E.D.
\smallskip
Our next theorem provides a lower bound for the number of generators in the
normal closure.

\medskip
\proclaim{\quad 2.8. Theorem}
For any $g \in G_S$ let $\tilde{\delta}(g) = \{x \in S \,|\, ord_x(g) \not\equiv 0 \,\roman{mod}\,2\}$.
Let $H$ be the normal closure in $G_S$ generated by a family of elements that
contains a subfamily of elements  $h= (h_1, \dots, h_i, \dots)$ such that the
following condition holds.
\smallskip
({\bf J}): For any $i$ there exist $x_i \in \tilde{\delta}(h_i)$ such that
$x_i \notin \tilde{\delta}(h_j)$ if $i \neq j$.
\smallskip
Then $H$ cannot be the normal closure in $G_S$ of less than $\roman{card}\,h$
generators.
\endproclaim

\smallskip
This theorem immediately implies

\medskip
\proclaim{\quad 2.9. Corollary}
In the situation of Theorem 2.8, $H$ cannot be the normal closure of a
finite number of elements if there is an infinite subsystem  $h$
satisfying ({\bf J}).
\endproclaim

\smallskip

{\bf Proof of Theorem 2.8.}
Define a map of $G_S$ into the vector space $\bold{F}_2^S$  as follows:
$$ 
\psi:\,  G_S \to V,\ \psi(g) = (\dots ,
ord_x(g)\,\roman{mod}\,2, \dots ). 
$$
\smallskip
It follows from 2.4(d) that  $\psi$ is a group homomorphism so that it maps conjugacy classes in $G_S$ into one element. The theorem will follow if
one shows that the image $\psi(H)$ cannot be generated  by less than
$\roman{card}\,h$ vectors.
But this follows immediately from the condition ({\bf J}) in the theorem
that guarantees that each image $\psi(h_i)$ has a non--zero
$x_{i}$--component  while all other vectors $\psi(h_j)$ have a zero
$x_{i}$--component. 

\medskip
\proclaim{\quad 2.10. Corollary}
None of the subgroups that are defined in (a),(c) and (d) in 2.7 can be
obtained as the normal closure of a finite number of generators.
\endproclaim

\smallskip
{\bf Proof.}
(a) follows from the fact that $\tilde{\delta}(t_x) = {x}$. Since $G_S$
contains the infinite number of $t_x$, $G_S$ cannot be obtained as the
normal closure of a finite number of elements.
\smallskip
(c) will follow similarly to (a) if we show that
$\tilde{\delta}(t_xt_yt_zt_{x^{\prime}}t_{y}t_{z'})$ contains  $x$ for $x \neq
y,z,x',z'$. This follows from the fact that the $t_x$--partition of
$t_xt_yt_z t_{x'}t_{y}t_{z'}$ is $(t_x, t_yt_z t_{x'}t_{y}t_{z'})$ (where
$t_yt_z t_{x'}t_{y}t_{z'} \in G_{S'}$).
Since $V(k)$ has infinitely many collinear triples $(x,y,z)$, such
that $x \neq y \neq z \neq x$, one can find infinitely many
generators in $B_0(V)$ satisfying the condition ({\bf J}).
\smallskip
The case (d) can be treated similarly.
\bigskip

\centerline{\bf \S3. Structure of universal equivalence}

\medskip

{\bf 3.1. Setup.} Let $P$ be an abstract cubic with  
the collinearity relation $L\subset P\times
P\times P$, such that for any $x,y\in P,$
there exists $z\in P$ with $(x,y,z)\in P.$

\smallskip

An equivalence relation $R$ on $P$ is called {\it admissible} if
the relation $L/R$ induced on $P/R$ has the following property:
for any $X,Y\in P/R,$ there exists a unique $Z$
with $(X,Y,Z)\in L/R.$ An admissible equivalence relation is called
{\it universal} if it is finer that any other
admissible relation.

\smallskip

In [M1] it was proved that the universal relation exists
(and of course, is unique) by a simple argument:
just take the intersection of all admissible relations.
Here we will clarify its structure by representing
it as a limit of a sequence of explicitly constructed
equivalence relations of which every next one is less
fine than the previous one.

\medskip

{\bf 3.2. Approximations.} For every $i\ge 0,$ we will
describe inductively a symmetric and reflexive binary
relation $\sim_i$ on $P$ and its transitive closure
$\approx_i.$ By definition, $\sim_0$ and $\approx_0$
are simply identical relations $x=x^{\prime}.$

\smallskip

{\bf 3.2.1. Definition.} {\it If $\sim_i$ and $\approx_i$
are already defined,
we put $x\sim_{i+1}x^{\prime}$ iff $x=x^{\prime}$ or there exist
$u,v,u^{\prime},v^{\prime}\in P$ such that
$u\approx_iu^{\prime},$ $v\approx_iv^{\prime},$
and $(u,v,x)\in L,$ $(u^{\prime},v^{\prime},x^{\prime})\in L.$

\smallskip

Furthermore, we put $x\approx_{i+1} x^{\prime}$ iff there is a sequence
of points $x=y_0, y_1,\dots ,y_r=x^{\prime}$ such that
$y_a\sim_{i+1} y_{a+1}$ for all $a<r.$}

\smallskip

Let us consider the case $i=1.$ By definition,
$x\sim_1 x^{\prime}$ iff there exist $u,v\in P$ such that
$(u,v,x), (u,v,x^{\prime})\in L.$ Let $P$ be the set
of $k$--points of a cubic surface $V$ and $L$ the usual
collinearity relation. Assume for simplicity
that $V$ does not contain lines defined over $k.$
Then $x\sim_1x^{\prime}$ means that
$x=x^{\prime}$ or
$x$ and $x^{\prime}$ lie on the intersection
of $V$ with the tangent plane at some $k$--point $u$
(with $u$ deleted if the double  tangent lines to $u$
in this plane are not defined over $k$). So one
equivalence class for $\approx_1$ consists of one point or
of a maximal connected union of such quasiprojective curves,
two of them being connected if they have
an intersection point defined over $k.$ The case of
general cubic surface allows a similar description, but
points of $k$--lines in $V$ must be added as subsets of
equivalence classes.

\medskip

{\bf 3.3. Theorem.}{\it (a) If $x\approx_i x^{\prime}$ then
$x\approx_{i+1} x^{\prime}.$

\smallskip

(b) Denote by $\approx$
the equivalence relation
$$
x\approx x^{\prime}\quad \Longleftrightarrow \quad
\exists i,\ x\approx_i x^{\prime}.
$$
Then it is admissible and universal.}

\smallskip

{\bf Proof.} (a) It suffices to prove that if $x\ne x^{\prime},$
$x\sim_i x^{\prime}$ then $x\approx_{i+1} x^{\prime}.$
For $i=0$ this is clear. Assume that we have proved
that $u\sim_{i-1} u^{\prime}$ implies $u\approx_i u^{\prime}.$

\smallskip

If $x\sim_i x^{\prime}$ then by definition
$(u,v,x)\in L$ and $(u^{\prime}, v^{\prime}, x^{\prime})\in L$
for some $u\approx_{i-1}u^{\prime}$ and $v\approx_{i-1}v^{\prime}.$
From the inductive assumption it follows that
$u\approx_{i}u^{\prime}$ and $v\approx_{i}v^{\prime}.$
By definition, then $x\approx_{i+1} x^{\prime}.$

\smallskip

(b) Let us first prove that $\approx$ is admissible,
in other words, if $(u,v,x)\in L,$
$(u^{\prime},v^{\prime},x^{\prime})\in L,$ and
$u\approx u^{\prime}, v\approx v^{\prime},$ then
$x\approx x^{\prime}.$ In fact,
for some $i$ we have $u\approx_i u^{\prime}, v\approx_i v^{\prime},$
so that $x\approx_{i+1} x^{\prime}$ and $x\approx x^{\prime}.$

\smallskip

Now denote temporarily the universal equivalence relation
by $\approx_U.$ The previous argument shows that
$x\approx_U x^{\prime} \Rightarrow x\approx x^{\prime}.$
It remains to prove that $x\approx_ix^{\prime} \Rightarrow
x\approx_U x^{\prime}.$ We argue by induction. Again, it suffices
to check that $x\sim_ix^{\prime} \Rightarrow
x\approx_U x^{\prime}$ assuming $x\ne x^{\prime}.$
We can then find $(u,v,x)\in L, (u^{\prime},v^{\prime},x^{\prime})\in L$
such that $u\approx_{i-1}u^{\prime},$ $v\approx_{i-1}v^{\prime}.$
Therefore $u\approx_U u^{\prime}$, $v\approx_U v^{\prime}$,
and finally $x\approx_U x^{\prime}$.

\medskip

{\bf 3.4. Types of finite generation.} Let us say, as in [M3], that $P$ is
$\circ$--generated by $(x_{\alpha}\,
|\,\alpha\in A)$ if for any $y\in P$ there is a non--associative
commutative word in
$x_{\alpha}$'s such that, informally, $y$ is one of the values
of this word. This means that when we calculate this word
in the order determined by the brackets, every time
that we have to calculate some $u\circ v,$ we may replace
it by any $x$ such that $(u,v,x)\in L.$

\medskip

{\bf 3.4.1. Claim.} {\it If $P$ is $\circ$--generated by $(x_{\alpha}
|\,\alpha\in A)$, then the CH--quasigroup $P/U$
is generated by the classes $X_{\alpha}$ of $x_{\alpha}.$}

\smallskip

We consider the following different types of $\circ$-generation.
\smallskip
{\bf 3.4.2. Values of nonassociative words.}
Let $W$ be a non-associative commutative word in finite number of variables
$X_i$, $P$ as in 3.1, and $x_i$ a family of elements of $P$
with the same set of indices.
We define different rules of computing values of $W$ on $(x_i)$ in the
order determined by the brackets inductively as follows for $i=0,1,
\dots\infty $. We set
$x \approx_{\infty} y$ if $x \approx y$ (i.e. $x \approx_j y$ for some
$j$).
\smallskip
{\bf Rule $A_i$.}
If the word $W=X$ has length 1, then a value of $W$ at any point $x
\in P$ is any $y \in P$ such that $x \approx_i y$. In
particular, $A_0$ means that the value of $W$ at $x$ coincides with $x$.
The rule $A_1$ means that the set of values of $W$ consists of those $y$ 
for which which there are points $u_j , y_j$, $j=0,\dots ,r$, $y_0 = x,
y_r = y$ such that  the following holds:
$(u_j, u_j, y_{j-1}) \in L, (u_j, u_j, y_j) \in L$ for $j=1,\dots ,r$.
\smallskip
If the word $W = X \circ Y$ has length 2, its set  $P(x,y)$ of values of
$W$ at $x,y \in P^2$ is defined as follows.
$$
P(x,y) = \{z \in P \,|\, z \approx_i z' , (x,y,z') \in L\}.
$$
If the word $W$ has  length more than two, it is a product of two non empty
words $W = W_1 \circ W_2$. Let $P(W_i)$ be a set of values of $W_i$ that is
defined inductively. Then the set of values $P(W)$
is defined as  $\cup P(x,y)$ for all
$(x,y) \in P(W_1) \times P(W_2)$.
\smallskip
We say that $P$ is $\circ_{A_{i}}$ generated by $P'=(x_\alpha \,|\, \alpha \in
A)$ if it is generated by application of  the rule $A_i$ to points in $P'$.
\smallskip
The inverse statement of 3.4.1 is valid for $\circ_{A_{\infty}}$ by trivial
reasons.
\medskip
{\bf 3.4.2. Claim.} {\it If CH-quasigroup $P/U\approx$ is $\circ$-generated
by classes of $(x_{\alpha}\,|\,\alpha\in A)$, then $P$ is
$\circ_{A_{\infty}}$ generated by  $x_{\alpha}.$}
\medskip
{\bf 3.4.3. Questions.}
 Let us define the {\it generation index } $i(P)$ of $P$ as the smallest $i$
such that  $P$ is $\circ_{A_i}$--generated by a finite number of points in $P$.
Let $P = V(k)$  for some cubic surface.
\smallskip
(1) For which fields $k$  and for which classes of cubic surfaces $i(P)$ is
finite? In particular, is $i(P)=0$ for $V$ defined over a number field
(the original Mordell-Weil problem)?
\smallskip
(2) If the the CH--quasigroup $P/U$ is finite, is the index $i(P)$ finite?
\smallskip
It would be worthwhile to study (2) for an abstract cubic $P$ that has an additional
property: every three points of it generate an Abelian group like points
on a plane cubic curve.
\bigskip

\centerline{\bf \S4. A group--theoretic description of universal
equivalence}

\medskip

In [M1], II.13.10 a group--theoretic description of universal equivalence
was given for a cubic surface that is defined over an infinite field and
has a point of general type.
In this section we extend this description of universal
equivalence. We relate the sequence of explicitly constructed equivalence
relations from \S 3 to a filtration by subgroups in the reflection
group associated with a minimal cubic surface.

\medskip
Let $B(V)$ and $B_0(V)$ be the groups  described in the examples 2.7.
Here the field $k$ over which the cubic surface $V$ is defined can be
finite and therefore we do not assume that $V(k)$ is infinite.

\smallskip
Define $x \sim y\,\roman{mod}\, U$ if $t_xt_y \in B_0(V)$. It is clear that $U$ is
an equivalence relation on $V(k)$. The proof of the following theorem differs
from the proof of the corresponding theorem 13.10 in [M1] in the following
respects. It uses the explicit description of the universal admissible
equivalence from the section 3 and the structural description of the
reflection group of $S=V(k)$.

\medskip

\proclaim{\quad 4.1. Theorem}
$U$ is the universal admissible equivalence relation.
\endproclaim

\smallskip
{\bf  Proof}. We will check in turn that each of the
equivalence relations is finer than the other one.

\smallskip

Assume first that $z'$ and $z$ are
universally equivalent. We want to show that
$z' \sim z\,\roman{mod}\, U$.

\smallskip
 According to Theorem 3.3,
$z' \approx_i z$ for some $i$.
Since $U$ is an equivalence relation,  it is sufficient to treat the case 
 $z' \sim_i z$. The following Lemmma does the job.

\medskip

\proclaim{\quad 4.2. Lemma} Denote by $B^i(V),
i=0,1,\dots ,$ the normal closure of the family $\{t_xt_{x'}\,|\,x \sim_i x'\}$
in $B(V)$. Let $x \sim_i x'$, $y \sim_i y'$, $(x,y,z) \in L$ and
$(x',y',z') \in L$ . Then the following holds:
$$
t_zt_{z'} \in t_zt_{z'}B^i(V) = t_xt_yt_zt_{x'}t_{y'}t_{z'}B^i(V)
\subset B^{i+1} \subset B_0(V)
$$
\endproclaim

\smallskip
{\bf Proof of Lemma 4.2.}
Using relations $t_x^2=1$ and $t_xt_yt_z = t_zt_yt_x$ we get $b =
t_zt_yt_xt_{x'}t_{y'}t_{z'} = t_{z}t_{z'}b'$ where $b' =
t_{z'}t_yt_xt_{x'}t_{y'}t_{z'}$. Next, $b'$ is conjugate to $b''=
t_yt_xt_{x'}t_{y'}$. And, finally, $b''$ is a product of $t_y t_{y'} \in
B^i(V)$ and $t_yt_xt_{x'}t_y$ which is conjugate to $t_xt_{x'} \in B^i(V)$.
This proves the equality $t_zt_{z'}B^i(V) =
t_xt_yt_zt_{x'}t_{y'}t_{z'}B^i(V)$. It remains
to show the inclusion $B^{i+1}(V) \subset B_0(V)$.
We will prove this inductively.

\smallskip

$B^1(V)$ is generated by $t_zt_{z'}$ such that $z'$ and $z$ lie on the intersection of $V$ with a tangent plane at some $k$--point $u$. In this
case $t_{z'}t_z = t_{z'}t_ut_ut_zt_ut_u \in B_0(V)$.

\smallskip

Assume that we already proved that $B^i(V) \subset B_0(V)$ and let us prove
that $t_zt_{z'} \in B_0(V)$. Let $z'' \in V(k)$ be such that $(x',y,z'')
\in L$.
Then $t_zt_{z'} = t_zt_{z''}t_{z''}t_{z'}$ and the following inclusions
hold:
$$ 
t_zt_{z''} \in t_zt_yt_xt_{x'}t_yt_{z''}B^i(V) \subset B_0(V) B^i(V)
\subset B_0(V),
$$
$$ 
t_{z''}t_{z'} \in t_{z''}t_{x'}t_yt_{z'}t_{x'}t_{y'}B^i(V) \subset
B_0(V) B^i(V) \subset B_0(V).
$$
Since $t_zt_{z'} \in B^{i+1}(V)$, this proves the inductive statement,
establishes the Lemma and the first part of the Theorem.

\smallskip
We turn now to the second part.
Let $A$ be any admissible equivalence relation. We shall show that
$x\sim y\,\roman{mod}\, U$ implies $x\sim y\,\roman{mod}\, A$.  Let $X,Y,Z$ be
the $A$--classes of
$x,y,z$. Then $Z = X
\circ Y$ in the sense of the composition law induced by collinearity
relation on $S=V(k)$. Denote by $E = V(k)/A$ the set of classes with the induced
structure of the symmetric quasigroup.
Let $t_X: E \rightarrow E$ be the map $t_X(Y) = X \circ Y$. The map $t_x\mapsto t_X$ extends to
an epimorphism of groups $\varphi: B(V) \rightarrow T(E)$. We will show
that its
kernel contains
$B_0(V)$. Therefore if $t_xt_y \in B_0(V)$
then $\varphi(t_xt_y) = t_Xt_Y =1$. This implies that $t_X=t_Y$ and that
$X=Y$.
To prove this property of $\varphi$  we need
to extend the Theorem 13.1 (ii),(iii) in [M1] to our case. Recall that
the Theorem 13.1 uses assumptions for cubic hypersurfaces that implies the
fact  that every equivalence class is dense in the Zariski topology. This
is not true any more in general in our case.
\smallskip
\proclaim{\quad 4.3. Lemma}
(a) $\varphi: B(V) \rightarrow T(E)$ is well defined and is an
epimorphism of groups.
\smallskip
(b) In $T(E)$ the following equality holds: $t_Xt_Yt_Z = t_{Y\circ Y}$.
\endproclaim

\smallskip
{\bf Proof.}
(a)  Our proof is based on the representation of elements in $B_0(V)$ as
``minimal'' words in the group $K^S$, the free product of groups $Z_2$
generated by symbols $T_x$, one for each point $x$ with the relations
$T_x^2=1$ (cf. [K1], 2.6 and \S 6). In order to construct the homomorphism
$B(V) \rightarrow T(E)$, we first define the action of $B(V)$ on $E$. Denote
by $T_{x_1}T_{x_2}
\dots T_{x_n}$ a minimal representation in $K^S$ of some
$s \in B(V)$. Choose $Y \in E$ and put $s(Y) = X_1 \circ (X_2 \circ \dots (X_n
\circ Y) \dots )$ where $X_i$ are classes of $x_i$ in $E$. 
\smallskip

One can show that
this definition does not depend on the choice of a minimal
representation of $s$ in $K^S$. This can be done inductively on the length
of minimal words in $K^S$. All minimal words of length one
representing the same element in $B(V)$ coincide. Let us assume that the
statement is proved for minimal words of the length $i-1$. Consider now
two different minimal words $w = T_1\dots T_i$, $w' = T'_1\dots T'_i$ of the
length $i$ representing $s \in B(V)$. (Minimal words representing the
same element have the same length). If $T_i = T'_i$ then the action of
$w$ (resp. $w'$)  on $E$ can be factored through the actions of $T_i$ and
$w_1 = T_1\dots T_{i-1}$ (resp. $w'_1 = T'_1\dots T'_{i-1})$ . Since $w_1$ and $w'_1$
represent the same element in $B(V)$ and have the length $i-1$, the statement
follows by the inductive assumption.

\smallskip

Otherwise, if $T_i \neq T'_i$, consider a $T_i$--partition of $w'$ (it is
defined in the same way as $t_x$--partition above): $(R_1, T_i, \dots R_{k-1}, T_i,
R_k)$. From [K1] it follows that $R_k =
T_{u_1}T_{v_1}T_{u_2}T_{v_2}\dots T_{u_r}T_{v_r}$ and $(u_j,v_j, u) \in L$ for
all $j=1,\dots r$ and $T_u = T_{i}$.  Moreover, if we replace
$T_{i} R_k$ in $w'$ with $R'_k T_{i}$ where  $R'_k =
T_{v_1}T_{u_1}T_{v_2}T_{u_2}\dots T_{v_r}T_{u_r}$, then we get a new word
$w''$ that is already a minimal representation of $s$.
Since $w''$ and $w$  both end with the same element $T_i = T_u$,
they act in the same way on $T(E)$. In order to prove that $w'$ and
$w''$ also act identically on $T(E)$ it is enough to check that $T_uR_k$
and $R'_kT_u$ act in the same way on $T(E)$. This can be shown using the
fact that $t_{u_j}t_{v_j}t_{u} = t_{u} t_{v_j} t_{u_j}$.

\smallskip
To complete (a) we need to show that for any two elements $s_1, s_2 \in
B(V)$ and  $Z \in E$ we have $s_1(s_2(Z)) = (s_1s_2)(Z)$. We will prove this
statement by induction  on the sum of lengths of  minimal
representation of $s_1$ and $s_2$. The statement is obvious if $s_1$ has
length 0. Assume now that $s_1$ has a minimal representation $w_1 =
T_{x_1}\dots T_{x_i}$, $ i\ge 1$, and $s_2$ has a minimal representation $w_2 =
T_{y_1} \dots T_{y_k}$. If $w = w_1w_2$ is the minimal representation of $s =
s_1s_2$ than the action of $s$ on $E$ is defined via the action of $w$ by
the rule $X_1 \circ (\dots X_i \circ (Y_1 \circ \dots (Y_k \circ Z) \dots )$ where
$X_i$ (resp. $Y_j$) are the  classes of $x_i$ (resp. $y_j$) and $Z \in E$.
Therefore $s_1(s_2(Z)) = (s_1s_2)(Z)$. Assume now that $w_1w_2$ is not
minimal.
\smallskip
Consider first the case when there exists such minimal representation of
$w_1, w_2$ that  $T_{x_i} = T_{y_1}$ (i.e. the last element in $w_1$
coincides with the first element in $w_2$). Let $s_1' \in B(V)$ be represented
by $w_1 = T_{x_1}\dots T_{x_{i-1}}$ and $s_2' \in B(V)$ be represented by
$w'_2 = T_{y_1} \dots T_{y_{k-1}}$. Then $s_1'(s_2'(Z)) = s_1(s_2(Z))$ and one
can apply  the inductive statement to $s_1'$ and $s_2'$.
\smallskip
Otherwise, let us assume that the word  $w_1w_2$ has the following 
$T_x$--partition;
$$
R_1 T_x R_2 \dots T_x R_l R_{l+1} T_x R_{l+2} T_x \dots T_x R_m 
$$ 
where
$R_1 T_x R_2 \dots T_x R_l$ (resp. $ R_{l+1} T_x R_{l+2} T_x \dots T_x R_m$) is
a minimal partition of $w_1$ (resp. $w_2$). Since $w_1w_2$ is not minimal,
$T_x$ can be chosen in such a way that $R_l R_{l+1} = T_{u_1}T_{v_1}T_{u_2}T_{v_2}
\dots T_{u_r}T_{v_r}$, where $(u_s, v_s, x) \in L$ for $s=1,\dots ,r$. As in the
case of minimal words above one can replace $T_x R_l R_{l+1}$  in $w_1w_2$
with 
$$ 
T_{v_1}T_{u_1}T_{v_2}T_{u_2} \dots T_{v_r}T_{u_r}T_x
$$  
and obtain a
new word $w'$ that has the same action on $E$ that  $w_1w_2$. Since $w'$
has two subsequent elements $T_x$, we can split it into a product of $w_1'$ that
ends with $T_x$ and  $w_2'$ that starts with $T_x$. This case was
already considered in this proof. 

\smallskip

(b) follows from properties of the group law on plane cubic curves.
This proves the Lemma 4.2.
\smallskip
To finish the proof of Theorem 4.1, we use the following identity:
$$ 
\varphi (t_xt_yt_zt_{x'}t_{y}t_{z'}) = t_Xt_Yt_{X\circ
Y}t_{X'}t_{Y}t_{X'\circ Y'} = t_{Y\circ Y}^2 = 1.
$$ 
Here $X,Y, \dots$ are the
classes of $x,y, \dots \mod A$. As a consequence, $B_0(V) \subset Ker\, 
 \varphi$,
proving the theorem.
\medskip
\proclaim{\quad 4.4. Corollary}
Let $V$ be a minimal cubic surface over a finite field with $q$ elements.
Then $B(V)/B_0(V) = Z_2$, except when all points of $V(k)$ are Eckardt
points. In the later case we have either $q=2$, $\roman{card}\,V(k) = 3$, or
$q=4$, $\roman{card}\,V(k) = 9$.
\endproclaim
\smallskip
{\bf Proof.} This
follows from the description of the universal equivalence for $V$ over finite
fileds in [Sw--D].
\medskip
{\bf 4.5. Remarks.}
(a) As it follows from the proof of  Theorem 4.1, it can be extended to an
abstract cubic for which every three points generate an abelian group,
in the same sense as for a plane cubic curve. We believe that this
theorem can be proved also for an abstract cubic using only a structural
description of $G_S$ without this additional assumption. We plan to address
this problem elsewhere.
\smallskip
(b) Groups $G_S$ were studied in [P] using different methods. [P] asked
whether {\it the dependency problem} $DP(n)$ is decidable for reflection groups of
an abstract cubic for $n\ge 3$ or $n= \infty$. $DP(n)$ can be formulated as
follows.
\smallskip
We will say that $g_0$ is {\it dependent} on $(g_1, \dots, g_k)$ if there is
a family $(g_{i_1}, \dots ,g_{i_p}$ and elements $u_1, \dots , u_p$ of $G$
such that
$$ 
g_0(u_1g_{i_1}u_1^{-1}) \dots (u_pg_{i_p}u_p^{-1}) = 1.
$$
If $n$ is a positive number or infinity then the {\it dependence problem}
$DP(n)$ asks for an algorithm to decide for any sequence $(g_0, \dots , g_k)$, $0\leq k < n$, of elements of $G$ whether or not $g_0$ is dependent on $(g_1,
\dots, g_k)$.
The problems $D(1), D(2)$ are usually called the {\it word problem} and the
{\it conjugacy problem}.
\smallskip
A special case of the dependence problem for $t_xt_y \in B_0(V)$  can be
related to the decidability of universal equivalence. Namely, if
$DP(\infty )$ is decidable for $g_i =
t_{x_i}t_{y_i}t_{z_i}t_{x'_i}t_{y_i}t_{z'_i}$
and $g_0 = t_xt_y$ than one can efficiently define whether $x,y$ are
universally equivalent.
\smallskip
Since the decidability of the universal equivalence seems to be a very difficult
problem in general, one can infer about the difficulty of the  $DP(\infty)$
for $B_0(V)$.
\smallskip
{\bf Question.}
Let an abstract cubic $S$ be decidable. Is $DP(n)$ decidable for arbitrary
$t_xt_y$ and generators of the subgroup $B_0(V)$  described in 2.7(c)?
\medskip

(c) Another construction of a filtration of the group of
birational automorphism of $V$ reflecting the structure of admissible
equivalences is given  in [M2]. One can apply the method from [M2] to
the classes of universal equivalence. One can show that there exist classes of
universal equivalence that are abstract cubics. One can consider universal
equivalence on the set of points of such a class (considered as the abstract
cubic). Applying this construction iteratively one can get a set of
abstract cubics that corresponds to a filtration of subgroups in reflection
groups. As in [M2] one can ask whether this sequence of subgroups
stabilizes and what is its intersection.

\bigskip
\centerline{\bf \S 5. Birationally trivial cubic surfaces: a finiteness
theorem}

\medskip

{\bf 5.1. Modified composition.} Let $V$ be a smooth cubic surface,
and $x,y\in V(k).$ Let $C\subset V$ be a curve on $V$ passing
through $x,y$, and $p:\,C\to \bold{P}^2$ an embedding
of $C$ into a projective plane such that $p(C)$ is again a cubic, and
$p(x)\circ p(y)$ is defined in $p(C).$ We assume that $C$ and $p$
are defined over $k.$
In this situation, following [M3], we will put
$$
x\circ_{(C,p)}y:=p^{-1}(p(x)\circ p(y)).
$$

\smallskip

{\it Example 1.} Choose $C=\roman{\ a\ plane\ section\ of\ }V$
containing $x,y$. If $p$ is the embedding of $C$ into the secant
plane, then $x\circ_{(C,p)}y=x\circ y$ in the standard notation.
Notice that the result does not depend on $C$ if $x\ne y.$
If $x=y,$ then the choice of $C$ determines a choice
of one or two tangent lines to $V$ at $x$ so that the multivaluedness
of $\circ$ is taken care of by the introduction of this new parameter.

\smallskip

{\it Example 2.} Assume now that $V$ admits a birational morphism
$p:\,V\to \bold{P}^2$ defined over $k$ (e.g., $V$ is split). We will
choose and fix $p$ once for all.
Then any plane section $C$ of $V$ not containing one of the
blown down lines as a component is embedded by $p$ into $\bold{P}^2$
as a cubic curve. Therefore we can apply to $(C,p)$ the previous
construction. Notice that this time $x\circ_{(C,p)}y$ depends
on $C$ even if $x\ne y.$

\smallskip

\proclaim{\quad 5.2. Theorem} Assume that $k$ is a finitely generated
field.
In the situation of Example 2,
the complement to the blown down lines in $V(k)$ is finitely
generated with respect to operations $\circ_{(C,p)}$ with the additional
restriction:
\smallskip
({\bf C}) the operation  $x \circ_{(C,p)}y$ is applied only to the different
previously constructed points.
\endproclaim

\smallskip
{\it Proof.}
This theorem was stated and proved in [M3] without the additional condition
({\bf C}). It uses the following auxiliary construction.
Choose a $k$--rational line $l\subset \bold{P}^2.$ Then $\Gamma
:=p^{-1}(l)$
is a twisted rational cubic in $V.$ The family of all such
cubics reflects properties of that of lines: a) any two different
points $a,b$ of $V(k)$ belong to a unique $\Gamma (a,b)$;\ b) any two
different
$\Gamma$'s either have one common $k$--point, or intersect
a common blown down line. The proof of this theorem is based on generation
of points by adding intersections of lines $l$ passing through pairs of
previously constructed points in a projective plane. This induces
generation of points on $V$ that are intersections of $p^{-1}(l)$. Analysis
of this proof in [M3] shows that it considers only {\it different} points
in pairs of previously constructed points hereby providing the statement of
the theorem with the condition ({\bf C}).
\bigskip
If one drops the condition ({\bf C}) one can prove the stronger statement.
\proclaim{\quad 5.3. Theorem}
 Let $ V $ be a smooth cubic surface over an arbitrary field $ k.$ Assume
that
$ V $ admits a birational morphism $ p: V \to \bold{P}^2 $. Then the
complement
$ P $ to all blown down lines in $ V(k) $ is generated by any single point
from
$ P$ (in the sense of the composition $\circ_{(C,p)}$).
\endproclaim
\smallskip
{\bf Proof. }
Let us choose a point $ x \in P $. The theorem will follow if we prove that
the set of points $x \circ_{(C,p)} x$ contains $P $ (here $ C $ runs through all
$k$--rational
plane sections of $ V $ passing through $x $).
Let us show that for any other point $y$ in $P$ there exists such
$C$  that $y = x \circ_{(C,p)} x $.
Indeed, following arguments of [M3],  for $ y \in P $
there exists a twisted cubic curve $G(x,y) := p^{-1}(l) $ where $ l $ is
the line
through $p(x),\,p(y)$ in $\bold{P}^2$. Let $l_1$ be the tangent line to $G(x,y)$ at
$x$.
Let a plane through points $x,y$ and $l_1$ cut a curve $C$ on $V$.
Then $l_1$ is a tangent line to $C$ at $x $, i.e. $G(x,y)$ is
tangent
to $C$ at $x$. Hence $l$ in $\bold{P}^2$ is
tangent to
$p(C)$ at $p(x)$. Since this line $l$ passes through $p(y)$, on $p(C)$ we
have
$p(y) = p(x) \circ p(x) .$ 
 This gives $y \in x \circ_{(C,p)} x $ proving the statement.
\medskip
One can apply this theorem to  the proof of the triviality of the 3--component
of the universal equivalence on $P=V(k)$. 3--component of the universal equivalence
can be defined as the finest admissible equivalence $U_3$ for which the
following condition holds:
\smallskip
{\it For any class} $X \in P/U_3, X \circ X = X$.
\smallskip
Simillarly one can define the 2--component of the universal equivalence as the
finest admissible equivalence for which the following condition holds:
\smallskip
{\it For any class $X \in P/U_2, X \circ X = O$ for some fixed class $O \in P$}.
\smallskip
It follows from [M1] that $U = U_3 \cap U_2$, where $U$ denotes the universal
equivalence.
\proclaim{\quad 5.4. Corollary}
 Let $ V $ be a smooth cubic surface over an arbitrary field $ k$. Assume
that
$ V $ admits a birational morphism $ p: V \rightarrow P^2 $. Then $U_3$ is
trivial on $V(k)$.
\endproclaim

The corollary can be deduced from the following two lemmas.
\proclaim{\quad 5.5. Lemma}
Let $C$ be a smooth plane cubic curve defined over a field $k$ such
that $C(k)$ is non--empty.
Let $p$ be  another plane embedding of  $C$ over $k$. Then
$$ 
x \circ_{(C,p)} y := p^{-1}(p(x) \circ p(y)) = t^{-1}((t(x) \circ
t(y))
$$
 where $t \in \roman{Bir}\,C$ is some birational automorphism of $C$ over $k$
which can be represented as a product of reflections of $C$
defined over $k$.
\endproclaim
\smallskip
{\bf Proof.}
The statement easily follows from the following fact: $p$ can be decomposed
into a product of reflections of $C$ over $k$ and a projective
isomorphism of $C$ and $p(C)$. Indeed, let us choose a point $0 \in C(k)$.
Isomorphism classes of invertible sheaves of  degree 3 are parametrized by the jacobian of  $C$ of
degree 3, say, $T$, and $T$ is a principal homogeneous space over $C$. This
means that $C(k)$ acts transitively on $T(k)$, i.e. any two sheaves $L_1,
L_2$ differ by a translation by a point $a \in C(k)$. Any translation
is a product of two reflections, whereas a projective
isomorphism preserves collinearity.

\smallskip
\proclaim{\quad 5.6. Lemma} In the same notation,
for any two points $x,y \in C(k)$ the following holds:
$$
t^{-1}(t(x)\circ t(y)) \sim x \circ y\,\roman{mod}\, U_3 .
$$
\endproclaim
\smallskip
{\bf Proof.}
Let $t = t_{x_1} \dots t_{x_n}$ where $x_i \in C(k)$. It is enough to check
the statement for $n=1$ since the general statement can be obtained by
induction. Let $t = t_z$. We have:
$t^{-1}(t(x)\circ t(y)) = t_z(t_z(x) \circ t_z(y)) = z \circ ((z \circ x)
\circ (z \circ y)) = z \circ ((z \circ z) \circ (x \circ y)) \sim z \circ (
z \circ (x \circ y)) \mod U_3 \sim x \circ y \mod U_3 $.
Here we used $z \circ z \sim z \mod U_3 $. Q.E.D.
\smallskip
We can now deduce the Corollary 5.4.
\smallskip
Fix some $x \in P,$ where $P$ is the complement to all blown down lines in
$V(k)$.
By the Theorem 5.3, any point  $ z \in P$ can be represented as
$x\circ_{(C,p)}x$.
Let $z = x\circ_{C,p}x$ for some $C$. If $C$ is singular then all points on
$C(k)$ are equivalent $\mod U_3$ (this is a general property of any
singular plane cubic curve that does not have a line as a component).
Otherwise, by lemmas 5.5 and 5.6
$$
z =  x\circ_{(C,p)}x \sim x\circ x\,\roman{mod}\, U_3 \sim x\,\roman{mod}\, U_3.
$$
\smallskip

{\bf 5.7. Elimination of $\circ_{(C,p)}$.} The use of the
modified operation $\circ_{(C,p)}$ is somewhat annoying, and
we would like to replace it by the standard composition
$\circ$. For example,
in the setup of the Theorem 5.2 for any three points $x,y,z$ on a plane smooth
section $C \subset V$ the following equality holds:
$$
(x \circ_{(C,p)} y) \circ_{(C,p)} z = (x \circ y) \circ z . 
$$
This naturally leads to the question whether one can obtain the traditional
Mordell--Weil statement for the composition $\circ$ using our 
finiteness results for $\circ_{(C,p)}$ and some tricks like
the formula above.

\smallskip
The remaining part of the paper is dedicated
to the description of our, not altogether successful, attempts to eliminate $\circ_{(C,p)}$.
We reformulate the finiteness theorem above  in terms
that do not use explicitly  compositions $ \circ_{(C,p)}$ and a
morphism $p$ of a cubic surface into a projective plane. We only use the
standard
operation $ \circ $ and implicitly use some intersections of planes with
lines
that belong to this cubic surface.

\smallskip

Before we can state a new statement we need to define a new kind of
operation
on a cubic surface that involves lines belonging to this cubic surface.

\smallskip

\proclaim{\quad 5.7.1. Definition }
 Let $ V $ be a smooth cubic surface over an arbitrary field $ k$. Let
$\Lambda = \{l_1, l_2,m\}$ be three (not necessary k-rational) lines
belonging to $ V $ and such
that the following properties hold:
\smallskip
({\bf A}) $l_1$ and $l_2$ are skew lines (i.e. they do
not have a common point) and $m$ intersects $l_1$ and $l_2$.
\smallskip

Given a triple of lines $ \Lambda $ satisfying ({\bf A}) and an arbitrary
plane $ T $ not containing lines in
$V$, let us define a new composition of points $u$, and $w$ on $ T \cap  V $
as follows:
\smallskip
({\bf B})  $u \circ_{(T,\Lambda )} w = ( x \circ y ) \circ [z \circ 
(u \circ w)] ,$
\smallskip
where  $x = l_1 \cap T $,  $y =  l_2 \cap T$ and $z = m \cap T$.
\endproclaim

\smallskip

Of course, the point $u \circ_{(T,\Lambda )} w $ is not necessarily
$k$--rational
even $u,w$, and $T$ are $k$--rational. But there is a special case when the
composition
$ \circ_{(T,\Lambda )} $ produces rational points (over $k$) when $u,w$,
and $ T$
are defined over $k$ (whereas lines in $ \Lambda $ are not necessarily defined
over $k$).
This case is described in the following statement that reformulates
the Theorem 5.2  in terms of the composition
$ \circ_{(T,\Lambda )} $.

\smallskip

\proclaim{\quad 5.7.2. Theorem}
Let $V$ be a smooth cubic surface. Assume that $V$ admits a birational
morphism to a projective plane defined over $k$.
Assume that $k$ is finitely generated field. Then there exists
a triplet of lines $  $ on $ V$ satisfying the property ({\bf A})
such that the following statement holds: the complement to the blown
down lines in $V(k)$ is finitely generated with respect to operations
$ \circ_{(T,\Lambda )}$ with the additional restriction:
\smallskip
({\bf D}) the operation  $x \circ_{(T,\Lambda )}y$ is applied only to 
different previously constructed points. 
 (Here $\Lambda$ is fixed and $T$ runs through some set of
k-rational planes).
\endproclaim

\smallskip

Similarly, one can
reformulate
Theorem 5.3 in terms of new operations.

\smallskip

\proclaim{\quad 5.7.3. Theorem }
Let $V$ be a smooth cubic surface over an arbitrary field $k$. Assume that
$V$ admits a birational
morphism to a projective plane defined over $k$. Then the complement
$ P $ to all blown down lines in $ V(k) $ is generated by any single point
from
$ P $ in the sense of  compositions $\circ_{(T,\Lambda)}$ for some
fixed triple of lines $ \Lambda $ in $V$.
\endproclaim

\smallskip

Below we will show how to replace  operations
$ \circ_{(C,p)} $ by  operations $\circ_{(T,\Lambda)}$.

\smallskip

\proclaim{\quad  5.7.4. Lemma}
Let $V$ be a smooth cubic surface defined over a field $k$ and $\bar{k}$ be
an algebraic closure of $k$. Let $ p: V \rightarrow \bold{P}^2 $ be a birational
morphism over $\bar{k}$. Then there exists a triplet of lines $\Lambda$
satisfying the property
({\bf A }) such that for any plane section $C$ of $V$ not containing one
of the blown down lines as a component and for any two points $u, w \in
V(\bar{k})$ lying
on $C$ the following holds:
\smallskip
$ u \circ_{(C,p)} w = u \circ_{(T,\Lambda)} w$
where $T$ is a plane that cuts the curve $C$ on $V$.
\endproclaim

\smallskip

\proclaim{\quad  5.7.5. Corollary }
Assume that the birational morphism $p$ in Lemma 5.7.4  is defined over
$k$.
Then a triplet $\Lambda$ can be chosen in such a way that the point
$u \circ_{(T,\Lambda)} w$ is $k$--rational if $u,w$ and the plane $T$ are
$k$--rational.
\endproclaim

\smallskip

The proof of  Lemma 5.7.4 is a consequence of the following claims which
might be of independent
interest.

\smallskip

\proclaim{\quad 5.7.6. Claim }
In the conditions of Lemma 5.7.4, let $x,y,u, w$ be some points on $C$.
Then the following equality holds:
$$
u \circ_{(C,p)} w = (x \circ y) \circ [(x \circ_{(C,p)} y)\circ (u \circ w
)]. 
$$
\endproclaim

\smallskip

In other words, if we know how to compute $z = x \circ_{(C,p)} y$ at least
for some two
points $x,y$ in $C$ then operation $\circ_{(C,p)}$ for all other points in
$C$
can be computed in terms of $\circ$ only.

\smallskip

\proclaim{\quad 5.7.7. Claim}
In the conditions of Lemma 5.7.4, let $\Lambda = \{l_1,l_2,m\}$
 be a triplet of lines satisfying ({\bf A}) and such that $p(m)$ is a line on
the plane $\bold{P}^2$, and $l_1, l_2$ are blown down lines.
Let $x = l_1 \cap T $,  $y =  l_2 \cap T$ and $z = m \cap T$, where the plane
$T$ cuts
a curve $C$ on $V$.
Then $z = x \circ_{(C,p)} y$.
\endproclaim

\smallskip

In other words, one can easily compute an operation $\circ_{(C,p)}$ for
intersection
of lines $l_1$ and $l_2$ with a plane $T$. The result of this composition
is an intersection of a third line $m$ with $T$ !

\smallskip

To show that the Lemma 5.7.4 follows from these claims, it is sufficient to
note
the following. By Claim 5.7.6, the operation $u \circ_{(C,p)} w $ can be
replaced by
$ (x \circ y) \circ [(x \circ_{(C,p)} y)\circ (u \circ w )] $ where $x,y$
are any points
on $C$. There exists a triplet of lines $\Lambda$ on $V$ satisfying ({\bf
A}), such that $p(m)$ is a line on the plane $\bold{P}^2$, and $l_1, l_2$ are the blown
down lines. By the Claim 5.7.7, $x,y$ can be chosen as intersections of lines
$l_1,l_2$ with a plane $T$ that cuts $C$ on $V$ and in this case $x
\circ_{(C,p)} y = m \cap T$.
\smallskip 
Now we prove our Claims.
\smallskip
{\bf Proof of the Claim 5.7.6.} {\it Step 1:}
Since $C$ and $p$ are fixed, one can simplify our notation by putting
$x * y =:x \circ_{(C,p)} y$. In this step  we show that for any points
$x,y,u,w$ on
$C$ the following equality holds:
$$ 
u*w = (x*y)(x \circ y)^{-1}(u \circ w),
\eqno(5.1)
$$
where the expressions in brackets
are multiplied by using an Abelian structure on $C$:
$ xy = a \circ (x \circ y) $
for some point $a$ in $C(k)$.

\smallskip
First, we consider the case when $C$ is smooth. In this case by the Lemma 5.5
 $p$ in the
formula $p(p^{-1}(u)\circ p^{-1}(w))$ can be replaced by a product of
reflections of $C$. Let us check  (5.1) for the case when $p$
can be replaced by one reflection $t_b$:
$$
u*w = p(p^{-1}(u)\circ p^{-1}(w)) = b \circ ((b \circ u ) \circ (b \circ
w)) =
b \circ ((b \circ b) \circ (u \circ w )). 
$$
The general case can be obtained by iterating this argument.

\smallskip
Using the identity
$u \circ w = (a \circ a) u^{-1} w^{-1}  $ 
we get:
$$
 u*w = b \circ ((b \circ b) \circ (u \circ w )) = b^{-1} (b \circ b) (u
\circ w ) 
$$
Similarly we have for other two points:
$x*y = b^{-1} (b \circ b) (x \circ y ) $.
Replacing $b^{-1} (b \circ b)$ with $  (x*y)(x \circ y )^{-1}$ in $b \circ
((b \circ b) \circ (u \circ w )) $  gives (5.1).
\smallskip
{\it Step 2:}
Replacing the Abelian multiplication operation in (5.1) by $a \circ
(\dots)$ we can rewrite (5.1) as as:
$$
u*w =  a \circ ( r \circ (u \circ w)), 
$$
 where $r = a \circ \{(x*y) \circ
[(a \circ a) \circ (x \circ y)]\} $.
Since the point $a$ is arbitrary, one can choose $a = x \circ y $. This gives
$r = x*y $ and immediately implies the formula in the Claim 5.7.6.
\smallskip
In order to complete the proof of the Claim we need to consider the case when
$C$ is a singular plane cubic curve that does not contain a line. This 
can be done by appealing to an obvious limiting construction
in the case of topological field $k$, or to a similar
argument using the Zariski
topology in general.

\medskip
{\bf Proof of the Claim 5.7.7.}
Since $l_1 , l_2 $ are blown down lines and $p(m)$ is a line in $\bold{P}^2,$ the points
$p(x),\, p(y),\, p(z)$ are intersections of the line $p(m)$ with the curve $p(C)$
in $\bold{P}^2$. This means that
on $p(C)$ we have $p(x) \circ p(y) = p(z)$. This is equivalent to the equality $z =
x \circ_{(C,p)} y$ in the Claim.

\bigskip

\centerline{\bf References}

\medskip

[K1] D.~S.~Kanevski. {\it Structure of groups, related to cubic surfaces},
Mat. Sb. 103:2, (1977), 292--308 (in Russian); English. transl. in Mat. USSR
Sbornik, Vol. 32:2 (1977), 252--264.

\smallskip

[K2] D.~S.~Kanevsky, {\it On cubic planes and groups connected with cubic
surfaces}. J. Algebra 80:2 (1983), 559--565.

\smallskip

[M1] Yu.~I.~Manin. {\it Cubic Forms: Algebra, Geometry, Arithmetic.} North
Holland, 1974 and 1986.

\smallskip

[M2] Yu.~I.~Manin. {\it On some groups related to cubic surfaces}. In:
Algebraic Geometry. Tata Press, Bombay, 1968, 255--263.

\smallskip

[M3] Yu.~I.~Manin. {\it Mordell--Weil problem for cubic surfaces}. 
In: Advances in the 
Mathematical Sciences---CRM's 25 Years (L. Vinet, ed.)  CRM Proc. and Lecture Notes, vol. 11, Amer.~Math.~Soc., 
Providence, RI, 1997, pp. 313--318.

\smallskip

[P] S.~J.~Pride. {\it Involutary presentations, with applications to
Coxeter groups, NEC-Groups, and groups of Kanevsky}. J. of Algebra 120 (1989),
200--223.

\smallskip

[Sw--D] H.~P.~F.~Swinnerton--Dyer. {\it Universal equivalence for cubic surfaces over
finite and local fields.} Symp.~Math., Bologna 24 (1981), 111--143.

\vskip2cm

{\it E-mail addresses:}

\smallskip
kanevsky{\@}us.ibm.com
\smallskip
manin{\@}mpim--bonn.mpg.de

\enddocument